\documentclass[12pt]{amsart}
\usepackage{amsfonts, amssymb, amsmath, amsthm}
\usepackage[mathscr]{eucal}
\usepackage[hmargin=1.2in, vmargin=1.5in]{geometry}
\newtheorem{theorem}{Theorem}[section]
\newtheorem{proposition}[theorem]{Proposition}

\newtheorem{corollary}[theorem]{Corollary}

\theoremstyle{definition}
\newtheorem{definition} [theorem]{Definition}

\newtheorem{remark}{Remark}

\numberwithin{equation}{section}
\usepackage{xcolor}
\newcommand{\va}{\mathbf a}

\newcommand{\vp}{\mathbf p}
\newcommand{\vq}{\mathbf q}
\newcommand{\vu}{\mathbf u}
\usepackage[normalem]{ulem}
\newcommand{\vx}{\mathbf x}

\newcommand{\vl}{\mathbf l}

\newcommand{\R}{\mathbb R}

\begin{document}

\title[Directional short-time special affine Fourier transform of distributions]{Directional short-time special affine Fourier transform of distributions}

\author[I. Berisha]{Ilir Berisha}
\address{Faculty of Civil Engineering, University of Prishtina, Agim Ramadani, Prishtina, 10000, Kosovo}
\email{ilir.berisha@uni-pr.edu}

\author[A. Ferizi]{Astrit Ferizi}
\address{Faculty of Mathematics and Natural Sciences, University of Prishtina, George Bush 31, Prishtina, 10000, Kosovo}
\email{ferizi.astrit@gmail.com}

\author[K. Hadzi-Velkova Saneva]{Katerina Hadzi-Velkova Saneva}
\address{Faculty of Electrical Engineering and Information Technologies, Ss. Cyril and Methodius University in Skopje, Rugjer Boshkovikj 18, Skopje, 1000, North Macedonia}
\email {saneva@feit.ukim.edu.mk}



\begin{abstract}

In this paper, we introduce and study the directional short-time special affine Fourier transform (DSTSAFT) in the multidimensional setting. Motivated by the need to combine localized time-frequency analysis with directional sensitivity and the flexibility of the special affine Fourier transform, the DSTSAFT provides a unified framework that encompasses established directional transforms as particular cases. We establish Parseval's identity and derive an inversion formula for the DSTSAFT. Furthermore, we prove that both the DSTSAFT and its synthesis operator are continuous mappings between the corresponding spaces of test functions. Building on these continuity results, we extend the DSTSAFT and its synthesis operator to the corresponding spaces of tempered distributions via duality. We also establish an important relationship between the DSTSAFT and the special affine Fourier transform (SAFT) and derive a desingularization formula. The resulting framework provides an analytical basis for studying directional and localized properties of tempered distributions within the SAFT setting.
\end{abstract}

\maketitle

\section{Introduction}

 Abe et al. \cite{AbeSheridan1994} introduced the six-parameter integral transform known as the special affine Fourier transform (SAFT). Their original motivation was rooted in optics, particularly in the study of operations on light-wave functions. The term “special affine” reflects the fact that the transform corresponds to a special affine transformation of the time-frequency coordinates, whose matrix component belongs to the special linear group $ SL(2,\mathbb{R})$, consisting of $(2\times2)$ real matrices with determinant one \cite{hasanbiswas}. The SAFT, also known as the offset linear canonical transform, is among the most general integral transforms in time-frequency analysis, unifying a broad class of classical transforms used in signal processing and optics \cite{NingLiu}. In particular, it generalizes the Fourier transform (FT), the fractional Fourier transform (FRFT) \cite{Namias1980, Almeida1997, Almeida, Ozaktas}, and the linear canonical transform (LCT) \cite{ShahTanary}. In \cite{cai}, Cai established the frequency domain formulation of the SAFT by deriving an explicit relationship between the input and output spectra and demonstrating the symmetry of the transform in both the spatial and frequency domains. By providing additional flexibility in representing shifts, modulations, and related affine operations, the SAFT extends beyond its classical special cases and has found applications in optical system analysis, signal processing, quantum mechanics, and digital communications \cite{cai, Abe1994, ZhiWeiZhang, BhandariZayed}.

Beyond its classical $L^2$-theory, the SAFT has also been extended to  {tempered distributions} {using} a duality approach \cite{Kumar1,kumar}. 
In \cite{Lakshmanan}, a class of Boehmians containing the space of tempered distributions was introduced, providing another setting for the study of the SAFT. These developments allow for further analysis of the properties of the SAFT in generalized function spaces. In particular, pseudo-differential operators associated with the SAFT were introduced in \cite{Kumar1,kumar}.

The multidimensional SAFT was proposed by Gao et al. in \cite{gao}. In that work, the convolution theorem was established, and sampling theorems for deterministic and random bandlimited signals were presented within the multidimensional SAFT domain. 
In \cite{NingLiu}, the characterization of {multivariate} dynamical sampling in a shift-invariant space was established via both the discrete and continuous SAFT. Moreover, by means of two illustrative examples, the authors showed that the $n$-dimensional SAFT provides superior signal recovery performance compared with the $n$-dimensional Fourier transform and the $n$-dimensional linear canonical transform. {Sampling and reconstruction in $L^p(\mathbb R^n)$ subspaces associated with the multidimensional SAFT were investigated in \cite{yang}.} This work includes the development of canonical convolution structures, dual basis functions, and iterative algorithms for nonuniform sampling.

The SAFT does not provide local time-frequency information. To overcome this limitation, Shah et al. \cite{shah} introduced the short-time special
affine Fourier transform (STSAFT), which provides simultaneous localization
in time and in the SAFT-frequency domain \cite{RuiLi}.
 A slightly different definition of the STSAFT was proposed by Lakshmanan et al. \cite{LakshmananRoopkumarZayed}. This new transform is constructed using a newly introduced SAFT convolution, in contrast to the windowed SAFT introduced in \cite{shah}, which is defined directly in terms of a translated and modulated analyzing window.

In higher dimensions, time-frequency analysis may require an additional component: besides time and frequency localization, directional information may also be included to capture the geometric structure of signals, especially those containing edges and anisotropic features. In \cite{Candes1999}, Candès introduced the ridgelet transform, which is obtained by applying the one-dimensional wavelet transform to the Radon transform. Inspired by this construction, Grafakos and Sansing \cite{GrafakosSansing} introduced a directional version of the STFT defined by composing the Radon transform with the one-dimensional short-time Fourier transform (STFT). A slightly different directional variant of the STFT, called the directional short-time Fourier transform (DSTFT), was introduced and studied by Giv in \cite{Giv2013}. A thorough distributional framework for the DSTFT was developed by Saneva et al. in \cite{SanevaAtanasova2016}. Recently, in \cite{ferizi1}, following a similar idea, the directional short-time fractional Fourier transform (DSTFRFT) was introduced, extending the FRFT-based analysis to the directional setting and developing its distributional theory, including a desingularization formula. The developed theory was further applied to the quasiasymptotic analysis of distributions in \cite{FeriziSaneva}.

Very recently, Brahim \cite{Brahim} introduced the directional short-time
free metaplectic transform associated with free symplectic matrices. 
Motivated by these developments, in this paper we introduce the
directional short-time special affine Fourier transform (DSTSAFT), which, unlike the zero-offset transform in \cite{Brahim}, incorporates the affine {offset parameters} of the SAFT. {Compared with existing directional short-time transforms, the DSTSAFT combines directional and localized analysis with the additional flexibility provided by the affine parameters of the special affine Fourier transform (SAFT), while retaining the DSTFT and the DSTFRFT as particular cases.}
To the best of our knowledge, this more general directional
short-time formulation of the SAFT, together with its distributional
framework, has not yet been systematically studied.
We establish its relation to the SAFT, prove Parseval's identity and an inversion formula, and show that the DSTSAFT and its synthesis operator act continuously on the corresponding Schwartz spaces. These continuity results enable the extension of both operators to
tempered distributions by duality. The distributional framework extends the applicability of the DSTSAFT beyond ordinary functions, allowing the treatment of  {tempered distributions}, including singular {ones.} In this setting, we also obtain a Calderón reproducing formula and a desingularization formula. 

This paper is organized as follows. Section \ref{Notations} collects the notation, the spaces of test functions and distributions, and the integral transforms that are used in our work. Then, in Section \ref{DSTSAFTsection}, we introduce the directional short-time special affine Fourier transform (DSTSAFT), establish  a representation of the DSTSAFT in terms of the SAFT, and prove  Parseval's identity and an inversion formula. In Section \ref{DSTSAFTofSchwartzspace}, we introduce the synthesis operator, establish the continuity of the DSTSAFT and its synthesis operator, extend both operators to appropriate spaces of distributions, provide an important relation between the duality definition and the SAFT, and conclude the section with a desingularization formula.
Finally, Section \ref{LS} provides the proofs of the main continuity theorems.

\section{Preliminaries}\label{Notations}

\subsection{Notation and spaces}

In this subsection, we recall some notations and spaces that are used in our work. For $\mathbf{x} = (x_1, x_2, \ldots, x_n)^{{T}}, \, \mathbf{y} = (y_1, y_2, \ldots, y_n)^{{T}} \in \mathbb{R}^n$ and a multi-index $ 
\boldsymbol{\alpha}= (\alpha_1, \alpha_2, \ldots, \alpha_n) \in \mathbb{N}_0^n$, we write $\mathbf{x}^{\boldsymbol{\alpha}} = x_1^{\alpha_1} x_2^{\alpha_2} \cdots x_n^{\alpha_n},\,
\boldsymbol{\alpha}! = \alpha_1! \alpha_2! \cdots \alpha_n!, \,
\partial_{\vx}^{\boldsymbol{\alpha}} = \frac{\partial^{|\alpha|}}{\partial x_1^{\alpha_1} \partial x_2^{\alpha_2} \cdots \partial x_n^{\alpha_n}}, \, 
|\boldsymbol{\alpha}| = \alpha_1 + \alpha_2 + \cdots + \alpha_n,$
\(|\mathbf{x}|\) denotes the Euclidean norm, and \(\mathbf{x} \cdot \mathbf{y}={\mathbf{x}^{T} \mathbf{y}} = x_1 y_1 + x_2 y_2 + \cdots + x_n y_n\) is the scalar product of \(\mathbf{x}\) and \(\mathbf{y}\). We write $k_1 \lesssim k_2$ when \(k_1 \leq k \, \cdot k_2\) for some constant \(k>0\).  The notation  {$(f, \varphi)_{L^2(\R^n)}=\int_{\mathbb{R}^n} f(\mathbf{x})\overline{\varphi(\mathbf{x})}d\mathbf{x}$} {stands} for the \(L^2(\R^n)\)-inner product of \(f\) and \(\varphi\), and \(\langle f, \varphi \rangle\) for the dual pairing between the distribution \(f\) and a test function \(\varphi\). {When it makes sense, we have $(f, \varphi)_{L^2(\R^n)} = \langle f, \overline{\varphi} \rangle.$}

The Fourier transform (FT) \(\mathcal{F}\) is given by
\begin{equation} \label{FT}
    \mathcal{F} f(\boldsymbol{\xi}) = \hat{f}(\boldsymbol{\xi}) = (2\pi)^{-n/2} \int_{\mathbb{R}^n} f(\mathbf{x})e^{-i \mathbf{x} \cdot \boldsymbol{\xi}} \, d\mathbf{x}, 
\quad  f \in L^1(\mathbb{R}^n).
\end{equation}

The Schwartz space \(\mathcal{S}(\mathbb{R}^n)\) consists of all functions \(\varphi \in C^{\infty}(\mathbb{R}^n)\) such that
\begin{equation*} \label{seminorm} 
\rho_{\nu}(\varphi) = \sup_{\substack{\mathbf{x} \in \mathbb{R}^n\\ |\boldsymbol{\alpha}| \leq \nu}} \left(1 + |\mathbf{x}|^2\right)^{\nu/2} |\partial_{\vx}^{\boldsymbol{\alpha}} \varphi(\mathbf{x})| < \infty, 
\end{equation*}
for every \(\nu \in \mathbb{N}_0\) \cite{hor, Schwartz}. Its topology is defined in the standard way. Its dual {space} \(\mathcal{S}'(\mathbb{R}^n)\) is the space of tempered distributions.

We denote by $\mathbb{S}^{n-1}$ the unit sphere in $\R^n$ and assume throughout that $n \ge 2$. We regard $\mathbb{S}^{n-1}$ as embedded in $\R^n$ and identify functions $\varphi$ on $\mathbb{S}^{n-1}$ with their homogeneous extensions $\varphi^{*}$ of degree $0$ to $\R^n \setminus \{0\}$. One may show that $\varphi$ is smooth on $\mathbb{S}^{n-1}$ if and only if $\varphi^{*}$ is smooth on $\R^n \setminus \{0\}$. In this sense, we consider  the Laplace-Beltrami operator $\Delta_{\mathbb{S}^{n-1}}$ {on $\mathbb{S}^{n-1}$}. The space $\mathcal{S}(\mathbb{Y}^{2n})$, where  $\mathbb{Y}^{2n} = \mathbb{S}^{n-1} \times \mathbb{R} \times \mathbb{R}^n$, consists of all smooth functions $\Phi \in C^\infty(\mathbb{Y}^{2n})$ such that 
\begin{equation} \label{norm}
\rho^{\vl,m,k}_{s,r}(\Phi) = 
\sup_{(\mathbf{u}, b, \mathbf{a}) \in \mathbb{Y}^{2n}}
\left(1 + |b|^2\right)^{r/2}\left(1 + |\mathbf{a}|^2\right)^{s/2}
\left| \partial_\mathbf{a}^{\vl} \partial_b^m \Delta_{\mathbb{S}^{n-1}}^k \Phi(\mathbf{u}, b, \mathbf{a}) \right| < \infty
\end{equation}
for all $s, r, m, k \in \mathbb{N}_0, \vl \in \mathbb{N}_0^n$. Since $\mathbb{S}^{n-1}$  is compact, decay conditions are imposed only in $b$ and $\mathbf{a}$, while smoothness in the directional variable $\mathbf{u}$ is controlled by iterates of the Laplace--Beltrami operator.
The topology of $\mathcal{S}(\mathbb{Y}^{2n})$ is defined in the usual way by the family of  seminorms \eqref{norm}, and its dual is denoted by $\mathcal{S}'(\mathbb{Y}^{2n})$. 

We equip $\mathbb{Y}^{2n}$  with the standard measure $db d\mathbf{a}d\mathbf{u} $, where  $d\mathbf{u}$ stands for the normalized surface area measure on the unit sphere $\mathbb{S}^{n-1}$, and $db$ and $d\va$ denote the Lebesgue measures on $\mathbb{R}$ and $\mathbb{R}^n$, respectively.
The inner product in the space $L^2\left(\mathbb{Y}^{2n}\right):=L^2\left(\mathbb{Y}^{2n},dbd\va d\vu\right)$ is given by {\cite{Giv2013, SanevaAtanasova2016}} 
\begin{equation*}
(F,G)_{L^2(\mathbb{Y}^{2n})}
=
\int_{\mathbb{S}^{n-1}}
\int_{\mathbb R^n}
\int_{\mathbb R}
F\left(\vu,b,\va\right)\overline{G\left(\vu,b,\va\right)}
dbd\va d\vu.
\end{equation*}
Every measurable function $F$ on $\mathbb{Y}^{2n}$ that satisfies the bound 
\begin{equation*}
    |F(\mathbf{u}, b, \mathbf{a})| \le C \left(1 + |\mathbf{a}|^2\right)^{s/2} \left(1 + |b|^2\right)^{r/2}, 
\enspace (\mathbf{u}, b, \mathbf{a}) \in \mathbb{Y}^{2n},
\end{equation*}for some $C>0, \, s,r\in \mathbb{N}_0$, defines an element of $\mathcal{S}'(\mathbb{Y}^{2n})$ via the duality 
\begin{equation}\label{duality}
\langle F, \Phi \rangle := \int_{\mathbb{S}^{n-1}} \int_{\mathbb{R}^n} \int_{\mathbb{R}}
F(\mathbf{u}, b, \mathbf{a}) \Phi(\mathbf{u}, b, \mathbf{a}) \, dbd\va d\vu,
\quad \Phi \in \mathcal{S}(\mathbb{Y}^{2n}).
\end{equation}
{The polynomial growth of $F$, together with the rapid decay of $\Phi$,
ensures that the integral in (\ref{duality}) is absolutely convergent and defines
a continuous linear functional on $\mathcal{S}(\mathbb{Y}^{2n})$.}

The Schwartz kernel theorem  (\cite{LSchwartz}, Prop. 28, page 98)  yields the following topological isomorphism
\begin{equation*}
    \mathcal{S}(\mathbb{Y}^{2n}) {\cong}\mathcal{S}(\mathbb{S}^{n-1} \times \mathbb{R}) \, \widehat{\otimes} \, \mathcal{S}(\mathbb{R}^n),
\end{equation*}
where \( X \, \widehat{\otimes} \, Y \) denotes the completion of the algebraic tensor product
\( X \otimes Y \) with respect to the projective \(\pi\)-topology or, equivalently, the injective  \(\varepsilon\)-topology \cite{treves, hertle}. {This equivalence follows from the nuclearity of the Schwartz spaces}  $\mathcal{S}(\mathbb{R}^n)$ and $\mathcal{S}(\mathbb{S}^{n-1} \times \mathbb{R})$. {Moreover,} (\cite{treves}, Prop. 50.7) yields the following isomorphisms
\begin{equation*}
    \mathcal{S}'(\mathbb{Y}^{2n}) \cong \mathcal{S}'\left(\mathbb{S}^{n-1} \times \mathbb{R}, \mathcal{S}'(\mathbb{R}^n)\right) 
\cong \mathcal{S}'\left(\mathbb{R}^n, \mathcal{S}'(\mathbb{S}^{n-1} \times \mathbb{R})\right).
\end{equation*}
The latter two spaces are spaces of vector-valued distributions, see, for instance \cite{treves}. These three spaces can be identified via the standard identification
\begin{equation*} \label{standardidentification} 
\left\langle F, \varphi \otimes \Psi \right\rangle 
= \left\langle \left\langle F, \varphi \right\rangle, \Psi \right\rangle
= \langle \langle F, \Psi \rangle, \varphi \rangle,
\quad 
\varphi \in \mathcal{S}(\mathbb{R}^n), \ 
\Psi \in \mathcal{S}(\mathbb{S}^{n-1} \times \mathbb{R}).
\end{equation*}

All dual spaces in the paper are equipped with the strong dual topology, i.e., the topology of bounded convergence \cite{treves}.

\subsection{The Radon transform {(RT)}}

The Radon {transform (RT)}  of {an integrable} function $f \in L^1(\mathbb{R}^n)$, denoted by $R f$, is defined by
\begin{equation*}
    R f_{\mathbf{u}}(p) = \int_{\mathbf{x} \cdot \mathbf{u} = p} f(\mathbf{x}) \, d\mathbf{x} 
= \int_{\mathbb{R}^n} f(\mathbf{x}) \delta(p - \mathbf{x} \cdot \mathbf{u}) \, d\mathbf{x}, \enspace \mathbf{u} \in \mathbb{S}^{n-1}, \,p \in \mathbb{R},
\end{equation*}
where $\delta$ is the Dirac-delta function \cite{hertle, ramm, helgason}. 
{Fubini's theorem implies that $R f \in L^1(\mathbb{S}^{n-1} \times \mathbb{R})$ for every 
$f \in L^1(\mathbb{R}^n)$.} 
The Fourier slice theorem establishes the following relation between the FT and the RT, stating that for every $f \in L^1(\mathbb{R}^n)$, 

\begin{equation*}
    \widehat{R f_{\mathbf{u}}}(p) = \left(2\pi\right)^{\frac{n-1}{2}} \, \widehat{f}(p \mathbf{u}),
\end{equation*}
where $\mathbf{u} \in \mathbb{S}^{n-1}$ and $p \in \mathbb{R}$ \cite{ramm}.
 
\subsection{{The SAFT and STSAFT}}

In this section, we provide the fundamentals of the multidimensional SAFT. For a thorough analysis, we refer to  \cite{gao,yang}.  

Let $M = {(A,B,\mathbf{p};C,D,\mathbf{q})=}
\begin{bmatrix}
A & B & \mathbf{p} \\
C & D & \mathbf{q}
\end{bmatrix}$, where \( A, B, C, D \) are \( n \times n \) real matrices, and \( \mathbf{p}, \mathbf{q} \) are \( n \times 1 \) column vectors. We assume that $M_1=\begin{bmatrix}
A & B\\
C & D
\end{bmatrix}$ is a symplectic matrix, i.e., \( M_1^T J M_1 = J \), where  $J = 
\begin{bmatrix}
0 & I_n \\
-I_n & 0
\end{bmatrix}$ and $I_n$ is the $n$-dimensional identity matrix. Throughout the paper, we suppose that \( \det B \neq 0 \).

The multidimensional SAFT of \( f \in L^1(\mathbb{R}^n) \) with respect to the matrix \( M \) is defined as
\begin{equation} \label{SAFT}
S_Mf(\mathbf{a}) = \int_{\mathbb{R}^n} f(\mathbf{x}) K_M(\mathbf{x}, \mathbf{a}) \, d\mathbf{x},
\end{equation}
\noindent where the kernel is given by
\begin{align}\label{kerneli}\nonumber
    K_M(\mathbf{x}, \mathbf{a})=& (2\pi)^{-n/2} |\det B|^{-1/2} \exp \bigg\{ 
\frac{i}{2} \big( 
\mathbf{a}^T D B^{-1} \mathbf{a}
+ \boldsymbol{\Omega} B^{-T} \mathbf{a}
\\&+ \mathbf{x}^T B^{-1} A \mathbf{x}
+ 2 \mathbf{p}^T B^{-1} \mathbf{x}
- 2 \mathbf{a}^T B^{-T} \mathbf{x}
\big) 
\bigg\}
\end{align}
for $\boldsymbol{\Omega} = 2 (B \mathbf{q} - D \mathbf{p})^T,
\,\mathbf{x} = (x_1, x_2, \cdots, x_n)^T\in \R^n, 
\,\mathbf{a} = (a_1, a_2, \cdots, a_n)^T\in \R^n$ \cite{gao, yang}. 

The multidimensional special affine Fourier transform (SAFT) provides a unified framework that includes several well-known transforms. For the matrix $M=(A,B,\mathbf{0};\\C,D,\mathbf{0})$ in \eqref{SAFT}, the transform reduces to the linear canonical transform (LCT) \cite{ShahTanary}. In particular, choosing $A=D=\operatorname{diag}(\cos\theta_1,\ldots,\cos\theta_n)$, $B=-C=\operatorname{diag}(\sin\theta_1,\ldots,\sin\theta_n)$, with $\theta_k \notin \pi \mathbb{Z}$ for $k\in \{1,2,...,n\}$, yields the fractional Fourier transform
(FRFT)  \cite{Almeida, Ozaktas}. The Fourier transform (FT) \eqref{FT} is a special case obtained by choosing $M=(0,I_n,\mathbf{0};-I_n,0,\mathbf{0})$.

If $f \in L^1(\mathbb{R}^n)$ such that $S_M f\in L^1(\mathbb{R}^n)$, then
\begin{equation} 
f(\mathbf{x})
= \int_{\mathbb{R}^n} S_Mf(\mathbf{a}) 
\, K_{M^{*}}(\mathbf{a}, \mathbf{x}) \, d\mathbf{a},
\label{eq:saft_inverse}
\end{equation}
{for a.e. $\mathbf{x} \in \mathbb{R}^n$} {(\cite{gao}, Rel. 2.3),}  where $K_{M^{*}}$ is the kernel corresponding to the  matrix
\begin{equation} \label{M*}
    M^{*}=
\begin{bmatrix}
D^{T} & -B^{T} & B\vq-D\vp \\[2mm]
-C^{T} & A^{T} & B^{-T}C^{T}B\vp-B^{-T}A^{T}B\vq
\end{bmatrix},
\end{equation}
and by (\cite{yang}, Lemma 1), $K_{M^{*}}(\va,\vx)=\overline{K_{M}(\vx,\va)}$ for $\vx,\va \in \R^n$. Furthermore, according to (\cite{gao}, Thm. 2.4), the following Parseval's identity for the SAFT holds
\begin{equation} \label{ParsevalSAFT}
\left(S_M f,\, S_M g \right) _{L^2(\mathbb{R}^n)} =\left( f,\, g \right) _{L^2(\mathbb{R}^n)},
\end{equation}
for all $f,g\in L^2(\R^n)$. By an argument analogous to the one-dimensional case (\cite{Kumar1}, Thm. 2.1), one can prove that the multidimensional SAFT is an isomorphism from $\mathcal{S}(\R^n)$ onto itself, and it may be extended to $\mathcal{S}'(\R^n)$ via the duality 
\begin{equation} \label{SAFTofD}
    \left\langle S_M f,\, \varphi \right\rangle
=\left\langle f,\, \overline{S_{M^{*}} \overline{\varphi}}\right\rangle,
\qquad f \in \mathcal{S}'(\mathbb{R}^n), \,\varphi \in \mathcal{S}(\mathbb{R}^n).
\end{equation}

\begin{remark}\label{Remark1}
The definition \eqref{SAFTofD} is consistent with the multidimensional SAFT for $f\in L^1(\R^n)$. {In contrast, the definition given in \cite[Def.~2.4]{Kumar1} is not consistent with the SAFT for integrable functions.} This inconsistency has also been pointed out in \cite{Lakshmanan}.
\end{remark}

In a natural way, and following the definition of the one-dimensional short-time special affine Fourier transform introduced in \cite{shah}, {we define} the multidimensional short-time special affine Fourier transform (STSAFT) of an integrable $f \in L^1(\mathbb{R}^n)$ with respect to a window $\psi \in L^{\infty}(\mathbb{R}^n)$ by 
\begin{equation}\label{STSAFT}
S_\psi^M f(\mathbf b,\mathbf{a})
=
\int_{\mathbb{R}^n}
f(\mathbf{x})\, \overline{\psi(\mathbf{x} - \mathbf{b})}\,
K_M(\mathbf{x},\mathbf{a})\, d\mathbf{x},
\end{equation}
where $\mathbf{b}, \va \in \mathbb{R}^n$, and $K_M(\mathbf{x},\mathbf{a})$ is given by (\ref{kerneli}).

The multidimensional STSAFT includes several well-known multidimensional time-frequency transforms as special cases. For the {parameter} matrix $M=(A,B,\mathbf{0};C,D,\mathbf{0})$ in \eqref{STSAFT}, we obtain the multidimensional short-time linear canonical transform (STLCT) \cite{ShahQadriLone}. In particular, the multidimensional short-time fractional Fourier transform (STFRFT) \cite{GaoLi} is obtained by choosing $A=D=\operatorname{diag}(\cos\theta_1,\ldots,\cos\theta_n)$, $B=-C=\operatorname{diag}(\sin\theta_1,\ldots,\\\sin\theta_n)$, with $\theta_k \notin \pi \mathbb{Z}$ for $k\in \{1,2,...,n\}$. Setting  $M=(0,I_n,\mathbf{0};-I_n,0,\mathbf{0})$ recovers the STFT \cite{Gröchenig}.

\section{The DSTSAFT} \label{DSTSAFTsection}

In this section, we introduce the DSTSAFT, establish an important relationship with the SAFT, and prove an extended Parseval's identity together with an inversion formula.
\begin{definition} \label{DefDSTSAFT} The directional short-time special affine Fourier transform (DSTSAFT) of $f \in L^1(\mathbb{R}^n)$ with respect to the window function $\psi{\in L^{\infty}(\R)}$ is defined by
\begin{equation} \label{DSTSAFT}
DS_{\psi}^M f(\mathbf{u}, b, \mathbf{a})
=
\int_{\mathbb{R}^n}
f(\mathbf{x})\, \overline{\psi(\mathbf{u} \cdot \mathbf{x} - b)}\,
K_M(\mathbf{x},\mathbf{a})\, d\mathbf{x},
\end{equation}
where $(\mathbf{u},b,\mathbf{a})\in \mathbb{Y}^{2n}$, and $K_M(\mathbf{x},\mathbf{a})$ is given by (\ref{kerneli}).
\end{definition} 
\begin{remark}
For the parameter matrix $M=(A,B,\mathbf{0};C,D,\mathbf{0})$ in \eqref{DSTSAFT} we obtain the definition of the directional {short-time free metaplectic transform introduced in \cite{Brahim}}. In particular, choosing $A=D=\operatorname{diag}(\cos\theta_1,\ldots,\cos\theta_n)$, $B=-C=\operatorname{diag}(\sin\theta_1,\ldots,\sin\theta_n)$, {with $\theta_k \notin \pi \mathbb{Z}$ for $k\in \{1,2,...,n\}$}, yields, up to a constant factor, the DSTFRFT introduced in \cite{ferizi1}. The special case $M=(0,I_n, \mathbf{0};-I_n, 0,  \mathbf{0})$ further recovers the DSTFT \cite{Giv2013}.    
\end{remark}

The following proposition establishes a relationship between the DSTSAFT and the SAFT.
\begin{proposition} \label{relationwithFT}
Let $\psi\in \mathcal S(\mathbb R)$ and $f\in L^1(\mathbb R^n)$. Then
\begin{equation*}
DS_\psi^M f(\mathbf{u},b,\mathbf{a})
=
(2\pi)^{-\frac12}
\int_{\mathbb R}
S_M\left(
f(\cdot)e^{-i(\cdot)\cdot(\xi \mathbf{u})}
\right)(\va)
\overline{\widehat{\psi}(\xi)}
e^{i\xi b}\,d\xi ,
\end{equation*}
where $(\mathbf{u},b,\mathbf{a})\in \mathbb{S}^{n-1}\times\mathbb R\times\mathbb R^n$.
\end{proposition}

\begin{proof} 
By using Fubini's theorem twice, the inversion formula for the FT, and the Fourier slice theorem, we obtain

\begin{align*}
DS_\psi^M f(\mathbf{u},b,\mathbf{a})
&=
\int_{\mathbb R^n}
f(\mathbf{x})K_M(\mathbf{x},\mathbf{a})
\overline{\psi(\mathbf{u}\cdot \mathbf{x}-b)}\,d\mathbf{x}
\\
&=
\int_{\mathbb R^n}
f(\mathbf{x})K_M(\mathbf{x},\mathbf{a})d\mathbf{x}
\int_{\mathbb R}
\overline{\psi(p-b)}
\delta(p-\mathbf{u}\cdot \mathbf{x})\,dp
\end{align*}\begin{align*}
&=
\int_{\mathbb R}
\overline{\psi(p-b)}dp
\int_{\mathbb R^n}
f(\mathbf{x})K_M(\mathbf{x},\mathbf{a})
\delta(p-\mathbf{u}\cdot \mathbf{x})\,d\mathbf{x}
\\
&=
\int_{\mathbb R}
R\left(f(\cdot)K_M(\cdot, \mathbf{a})\right)(\mathbf{u},p)
\overline{\psi(p-b)}\,dp
\\
&=
(2\pi)^{-\frac12}
\int_{\mathbb R}
R\left(f(\cdot)K_M(\cdot, \mathbf{a})\right)(\mathbf{u},p) dp
\int_{\mathbb R}
\overline{\widehat{\psi}(\xi)}
e^{i\xi (b-p)} d\xi
\\
&=
(2\pi)^{-\frac12}
\int_{\mathbb R}
\overline{\widehat{\psi}(\xi)}
e^{i\xi b}d\xi
\int_{\mathbb R}
R\left(f(\cdot)K_M(\cdot, \mathbf{a})\right)(\mathbf{u},p)
e^{-i\xi p}\,dp
\\
&=
(2\pi)^{\frac{n-1}{2}}
\int_{\mathbb R}
\mathcal F\left(f(\cdot)K_M(\cdot, \mathbf{a})\right)(\xi \mathbf{u})
\overline{\widehat{\psi}(\xi)}
e^{i\xi b}\,d\xi
\\
&=
(2\pi)^{-\frac12}
\int_{\mathbb R}
S_M\left(
f(\cdot)e^{-i(\cdot)\cdot(\xi \mathbf{u})}
\right)(\mathbf{a})
\overline{\widehat{\psi}(\xi)}
e^{i\xi b}\,d\xi .
\end{align*}

\end{proof}

The following result provides Parseval's identity for the DSTSAFT.

\begin{proposition}[Parseval's identity] \label{Parsevalidentity}
Let $\psi \in \mathcal{S}(\mathbb{R})\setminus\{0\}$. Then, for all $f,g\in L^1(\mathbb{R}^n)\cap L^2(\mathbb{R}^n)$, the following Parseval's identity holds
\begin{equation*}
    \int_{\mathbb{R}^n}
f(\mathbf{x})\overline{g(\mathbf{x})}\,d\mathbf{x}
=
\frac{1}{(\eta,\psi)_{L^2(\R)}}
\int_{\mathbb{S}^{n-1}}
\int_{\mathbb{R}^n}
\int_{\mathbb{R}}
DS_{\psi}^{M}f(\mathbf{u},b,\mathbf{a})
\overline{DS_{\eta}^{M}g(\mathbf{u},b,\mathbf{a})}
\,dbd\va d\vu,
\end{equation*}
where $\eta\in \mathcal{S}(\mathbb{R})$ is a synthesis window for $\psi$, namely, $(\eta,\psi)_{L^2(\R)}\ne0$.
\end{proposition}

\begin{proof}

By using Prop. \ref{relationwithFT}, {Parseval's identity for the one-dimensional Fourier
transform,} Fubini's theorem, and Parseval's identity for the SAFT {\eqref{ParsevalSAFT}}, we obtain
\begin{align*}
&\int_{\mathbb{S}^{n-1}}
\int_{\mathbb{R}^n}
\int_{\mathbb{R}}
DS_{\psi}^{M}f(\mathbf{u},b,\mathbf{a})
\overline{DS_{\eta}^{M}g(\mathbf{u},b,\mathbf{a})}
\,dbd\va d\vu
\\
=&
\int_{\mathbb{S}^{n-1}}
\int_{\mathbb R^n}
\int_{\mathbb R}
S_M\!\left(f(\cdot)e^{-i(\cdot)(\xi \mathbf{u})}\right)(\mathbf{a})
\overline{
S_M\!\left(g(\cdot)e^{-i(\cdot)(\xi \mathbf{u})}\right)(\mathbf{a})}
\,\widehat{\eta}(\xi) \overline{\widehat{\psi}(\xi)}
\,d\xi\,d\mathbf{a}\,d\mathbf{u}
\\
=&
\int_{\mathbb{S}^{n-1}} 
\int_{\mathbb R}
\overline{\widehat{\psi}(\xi)}
\,\widehat{\eta}(\xi)\,d\xi \, d\mathbf{u}
\int_{\mathbb R^n}
S_M\!\left(f(\cdot)e^{-i(\cdot)(\xi \mathbf{u})}\right)(\mathbf{a})
\overline{
S_M\!\left(g(\cdot)e^{-i(\cdot)(\xi \mathbf{u})}\right)(\mathbf{a})}
\,d\mathbf{a}
\\
=&
\int_{\mathbb{S}^{n-1}} \int_{\mathbb R}
\overline{\widehat{\psi}(\xi)}
\,\widehat{\eta}(\xi)\,d\xi\, d\mathbf{u}
\int_{\mathbb R^n}
f(\mathbf{x})e^{-i\mathbf{x}\cdot(\xi \mathbf{u})}
\overline{
g(\mathbf{x})e^{-i\mathbf{x}\cdot(\xi \mathbf{u})}
}
\,d\mathbf{x}
=
(\eta,\psi)_{L^2(\R)}(f,g)_{L^2(\R^n)}.
\end{align*}
\end{proof}

{The Parseval's identity for the DSTSAFT may be expressed as follows}
$$\left(
DS_\psi^M f,
DS_\eta^M g
\right)_{L^2(\mathbb{Y}^{2n})}=(\eta,\psi)_{{L^2(\R)}}(f,g)_{L^2(\mathbb R^n)}.$$
Moreover, if $\psi\in\mathcal{S}(\mathbb{R})$ satisfies $\|\psi\|_{L^2(\mathbb{R})}=1$ and $f \in \mathcal{S}(\mathbb{R}^n)$, then $\|DS_\psi^M f\|_{L^2(\mathbb{Y}^{2n})}
=
\|f\|_{L^2(\mathbb{R}^n)}$. Since $\mathcal{S}(\mathbb{R}^n)$ is a dense subset of $L^2(\mathbb{R}^n)$,  $DS_\psi^M$ 
can be extended to an isometric embedding $L^2(\mathbb{R}^n)\to L^2(\mathbb{Y}^{2n})$.\\

The next proposition gives the inversion formula for the DSTSAFT.

\begin{proposition}[The inversion formula] \label{Inversionformula}
Let $\psi\in\mathcal S(\mathbb R)$ be a non-trivial window and let
$\eta\in\mathcal S(\mathbb R)$ be a synthesis window for $\psi$.
Let $f\in L^1(\mathbb R^n)$ be such that $S_M f\in L^1(\mathbb R^n)$.
Then, for a.e. $\vx \in\mathbb R^n$,
\begin{equation*}
 f(\mathbf{x})=\frac{1}{(\eta,\psi)_{L^2(\R)}}
\int_{\mathbb{S}^{n-1}}
\int_{\mathbb R^n}
\int_{\mathbb R}
DS_\psi^M f(\mathbf{u},b,\mathbf{a})\,
\eta(\mathbf{x}\cdot \mathbf{u}-b)\,
K_{M^{*}}(\mathbf{a},\mathbf{x})
\,dbd\va d\vu.
\end{equation*}
\end{proposition}
\begin{proof} A direct calculation from (\ref{SAFT}) shows that
$S_M\left(f(\cdot)e^{-i(\cdot)\cdot\omega}\right)$ is, up to a unimodular factor,
a translate of $S_Mf$ by $B\omega$. Hence,
$S_M\left(f(\cdot)e^{-i(\cdot)\cdot\omega}\right)\in L^1(\mathbb{R}^n)$ whenever
$S_Mf\in L^1(\mathbb{R}^n)$, so that the SAFT inversion formula (\ref{eq:saft_inverse})
is applicable to $f(\cdot)e^{-i(\cdot)\cdot\omega}$. Therefore, by  using Prop. \ref{relationwithFT}, Fubini's theorem, and the inversion formula for the SAFT (\ref{eq:saft_inverse}), we obtain 
\begin{align*}
&\int_{\mathbb{S}^{n-1}}
\int_{\mathbb R^n}
\int_{\mathbb R}
DS_\psi^M f(\mathbf{u},b,\mathbf{a})\,
\eta(\mathbf{x}\cdot \mathbf{u}-b)\,
K_{M^{*}}(\mathbf{a},\mathbf{x})
\,dbd\va d\vu
\\
=&(2\pi)^{-\frac12}\int_{\mathbb{S}^{n-1}}
\int_{\mathbb R^n}
\int_{\mathbb R}
\eta(\mathbf{x}\cdot \mathbf{u}-b)\,
K_{M^{*}}(\mathbf{a},\mathbf{x}) 
\,dbd\va d\vu 
\int_{\mathbb R}
S_M\left(
f(\cdot)e^{-i(\cdot)\cdot(\xi \mathbf{u})}
\right)(\mathbf{a})
\overline{\mathcal F\psi(\xi)}
e^{i\xi b}\,d\xi
\\
=&
(2\pi)^{-\frac12}
\int_{\mathbb{S}^{n-1}}
\int_{\mathbb R}
\overline{\widehat{\psi}(\xi)} d\xi\,d\mathbf{u}
\int_{\mathbb R^n}
S_M\!\left(f(\cdot)e^{-i(\cdot)(\xi \mathbf{u})}\right)(\mathbf{a})
K_{M^{*}}({\mathbf{a},\mathbf{x}})
\,d\mathbf{a}
\int_{\mathbb R}
\eta(\mathbf{x}\cdot \mathbf{u}-b)e^{i\xi b}\,db
\\
=&f(\mathbf{x})
\int_{\mathbb{S}^{n-1}}
\int_{\mathbb R}
\overline{\widehat{\psi}(\xi)}
e^{-i\mathbf{x}\cdot(\xi \mathbf{u})}
e^{i\xi(\mathbf{u}\cdot \mathbf{x})}
\widehat{\eta}(\xi)
\,d\xi\,d\mathbf{u}\\=&
f(\mathbf{x})
\int_{\mathbb{S}^{n-1}}
\int_{\mathbb R}
\overline{\widehat{\psi}(\xi)}
\widehat{\eta}(\xi)
\,d\xi\,d\mathbf{u}=
f(\mathbf{x})(\eta,\psi)_{{L^2(\R)}}.
\end{align*}
\end{proof}

\begin{remark}
Although the results of this paper are developed under the assumption $n \geq2$, it is worthwhile to consider the previous proposition for the case $n=1$. In this case, $\mathbb{S}^{0}=\{1,-1\}$, and the integral over $\mathbb{S}^{0}=\{1,-1\}$ is equal to one-half of the sum of two integrals corresponding to {$u=1$ and $u=-1$}. Therefore, 
\begin{align*}
&\int_{\mathbb{S}^0}\int_{\mathbb R}
\int_{\mathbb R}
DS_{\psi}^{M}f(u,b,a)\eta(x\cdot u-b)\,
K_{M^{*}}(a,x)\,
dbda{du} \\
=&\frac{1}{2} \int_{\mathbb R}
\int_{\mathbb R}\left(
DS_{\psi}^{M}f(1,b,a)\eta(x-b)+DS_{\psi}^{M}f(-1,b,a)\eta(-x-b)\right)\,
K_{M^{*}}(a,x)\,
dbda. \end{align*}
Making the change of variable $b\mapsto -b$ in the term corresponding to $u=-1$, we obtain

\begin{align*}
&\int_{\mathbb{S}^0}\int_{\mathbb R}
\int_{\mathbb R}
DS_{\psi}^{M}f(u,b,a)\eta(x\cdot u-b)\,
K_{M^{*}}(a,x)\,
dbda{du} \\
=&\frac{1}{2} \int_{\mathbb R}
\int_{\mathbb R}
S_{\psi}^{M}f(b,a)\eta(x-b)K_{M^{*}}(a,x)dbda+\frac{1}{2}\int_{\mathbb R}
\int_{\mathbb R} S_{\check{\psi}}^{M}f(b,a)\check{\eta}(x-b)\,
K_{M^{*}}(a,x)
dbda,
\end{align*}
where $\check{\psi}(x)=\psi(-x)$ and $\check{\eta}(x)=\eta(-x)$. {Since $(\check\eta,\check\psi)_{L^2(\R)}=(\eta,\psi)_{L^2(\R)}$,} the STSAFT inversion formula (\cite{shah}, Thm. 2.6), applied to both terms,  yields 
\begin{equation} \label{15}
    \begin{split}
&\int_{\mathbb{S}^0}\int_{\mathbb R}
\int_{\mathbb R}
DS_{\psi}^{M}f(u,b,a)\eta(x\cdot u-b)\,
K_{M^{*}}(a,x)\,
dbda{du}\\=& \int_{\mathbb R}
\int_{\mathbb R}
S_{\psi}^{M}f(b,a)\eta(x-b)K_{M^{*}}(a,x)dbda=(\eta,\psi)_{L^2(\R)}f(x) 
\end{split}
\end{equation}
for a.e. $x\in \R$. The relation \eqref{15} establishes that, in the one-dimensional case, the inversion formula for the DSTSAFT coincides with that for the STSAFT. 
\end{remark}


\section{The DSTSAFT on Schwartz spaces of test functions and distributions}\label{DSTSAFTofSchwartzspace}

This section focuses on the continuity properties of the DSTSAFT on the Schwartz space of test functions and its extension to tempered distributions. The proofs of continuity results, specifically Theorems \ref{ContinuityDSTSAFT} and \ref{ContinuityDSTSAFTsynthesis}, will be addressed in Section \ref{LS}.

\subsection{Continuity theorems}\label{ContDSTSAFT}

We now establish the continuity of the DSTSAFT on $\mathcal S(\mathbb R^n)$.

\begin{theorem} \label{ContinuityDSTSAFT}Let $M = {(A,B,\mathbf{p};C,D,\mathbf{q})}$. The operator $DS^M:\mathcal S(\mathbb R^n)\times \mathcal S(\mathbb R)\rightarrow \mathcal S(\mathbb{Y}^{2n})$ defined by $(f,\psi) \mapsto DS_\psi^M f$ is a continuous sesquilinear mapping.
\end{theorem}

We aim to construct a left inverse of the DSTSAFT. For this purpose, we define the directional short-time special affine synthesis operator of  $\Phi\in \mathcal{S}(\mathbb{Y}^{2n})$ with respect to $\psi\in L^{\infty}(\R)$ as
\begin{equation}\label{synthesis}
\left(DS_\psi^M\right)^{*}\Phi(\mathbf{x})
:=
\int_{\mathbb{S}^{n-1}}
\int_{\mathbb R^n}
\int_{\mathbb R}
\Phi(\mathbf{u},b, \mathbf{a})
\psi(\mathbf{u}\cdot \mathbf{x}-b)
K_{M^{*}}(\mathbf{a},\mathbf{x})
\,db d\va d\vu,
\quad \mathbf{x}\in\mathbb R^n.
\end{equation}
It can be observed that the last integral is absolutely convergent for every $\vx\in\R^n$. Our next objective is to establish the continuity of the directional short-time special affine Fourier synthesis operator.

\begin{theorem}\label{ContinuityDSTSAFTsynthesis} Let $M = {(A,B,\mathbf{p};C,D,\mathbf{q})}$. The operator
$\left(DS^{M}\right)^{*}:
\mathcal S(\mathbb Y^{2n})\times\mathcal S(\mathbb R) \rightarrow \mathcal S(\mathbb R^{n}) $ defined by $(\Phi,\psi) \mapsto\left(DS_{\psi} ^{M}\right)^{*}\Phi$ is a  continuous bilinear mapping.
\end{theorem}

The following inversion formula is an immediate consequence of Prop. \ref{Inversionformula} and Theorems \ref{ContinuityDSTSAFT} and \ref{ContinuityDSTSAFTsynthesis}.

\begin{proposition}\label{Inversionsynthesis}
Let $\eta\in\mathcal S(\mathbb R)$ be a synthesis window for $\psi\in\mathcal S(\mathbb R) \setminus \{0\}$. Then, the following inversion formula holds for every $f\in \mathcal{S}(\R^n)$
\begin{equation}\label{Inversionsynthesisnr}
\left(\left(DS_\eta^{M}\right)^{*}\circ
DS_\psi^{M}\right)f(\vx)
= (\eta,\psi)_{L^2(\mathbb R)} 
f(\vx), \quad \vx \in \R^n.
\end{equation}
\end{proposition}
\begin{proof}
Since $f\in\mathcal S(\mathbb R^n)$ and $\psi\in\mathcal S(\mathbb R) \setminus \{0\}$, Thm. \ref{ContinuityDSTSAFT} implies that $DS_\psi^{M}f\in\mathcal S(\mathbb Y^{2n})$. {Moreover, since the SAFT is an isomorphism of
$\mathcal S(\mathbb R^n)$ onto itself,
$f,S_Mf\in\mathcal S(\mathbb R^n)\subset L^1(\mathbb R^n)$.
Hence, the assumptions of Prop.~\ref{Inversionformula} are satisfied.} Therefore, by the definition of the synthesis operator \eqref{synthesis} and Prop. \ref{Inversionformula}, we obtain, for every $\vx\in \R^n$,
\begin{align*}
\left(DS_\eta^{M}\right)^{*}DS_\psi^{M}f (\vx) 
&=\int_{S^{n-1}}
\int_{\mathbb R^n}
\int_{\mathbb R}
DS_\psi^{M}f(\mathbf u,b,\mathbf a)\,
\eta(\mathbf x\cdot\mathbf u-b)\,
K_{M^{*}}(\mathbf a,\mathbf x)
\,db d\va d\vu\\
&= (\eta,\psi)_{L^2(\mathbb R)} f(\mathbf{x}).
\end{align*}
\end{proof}


\subsection{The DSTSAFT of tempered distributions}\label{DSTSAFTinS'}

We now extend the definitions of the DSTSAFT and its corresponding synthesis operator to appropriate spaces of distributions.

We define the DSTSAFT of $f\in\mathcal S'(\mathbb R^n)$ with respect to a window function $\psi\in\mathcal S(\mathbb R)$, denoted by $DS_\psi^{M}f$, via the duality 
\begin{equation} \label{D1}
    \left\langle
DS_\psi^{M}f,\Phi
\right\rangle
=
\left\langle
f,
\overline{\left(DS_\psi^{M}\right)^{*}\overline{\Phi}}
\right\rangle,
\qquad
\Phi\in\mathcal S(\mathbb Y^{2n}).
\end{equation}

By Thm.~\ref{ContinuityDSTSAFTsynthesis}, the distribution $DS_{\psi}^{M}f$ is well defined via the duality \eqref{D1}. Furthermore, by (\cite{treves}, Corollary of Prop. 19.5), the mapping
$
DS_{\psi}^{M}:\mathcal{S}'(\mathbb{R}^n)\longrightarrow\mathcal{S}'(\mathbb{Y}^{2n})
$
is continuous. \\

The following result shows that, for $f\in L^1(\mathbb R^n)$, Def. \ref{DefDSTSAFT} and the definition given by \eqref{D1} coincide.

\begin{proposition}\label{relationsynthesis}
Let $\psi\in\mathcal S(\mathbb R)$, $f\in L^1(\mathbb R^n)$, and
$\Phi\in\mathcal S(\mathbb{Y}^{2n})$. Then
\begin{equation*}
    \left\langle
f,
\overline{\left(DS_\psi^{M}\right)^{*}\overline{\Phi}}
\right\rangle
=
\int_{\mathbb{S}^{n-1}}
\int_{\mathbb R^n}
\int_{\mathbb R}
DS_\psi^M f(\mathbf{u},b,\mathbf{a})\,
\Phi(\mathbf{u},b,\mathbf{a})
\,db d\va d\vu.
\end{equation*}
\end{proposition}

\begin{proof}
By Def. \ref{DefDSTSAFT}, relation \eqref{synthesis}, and Fubini's theorem, we have 
\begin{align*}
&\left\langle
f,
\overline{\left(DS_\psi^{M}\right)^{*}\overline{\Phi}}
\right\rangle=\int_{\mathbb R^n}
f(\mathbf{x})\overline{\left(DS_\psi^M\right)^*\overline{\Phi}(\mathbf{x})}\,d\mathbf{x}
\\
=&
\int_{\mathbb R^n}
f(\mathbf{x})d\mathbf{x}
\overline{
\int_{\mathbb{S}^{n-1}}
\int_{\mathbb R^n}
\int_{\mathbb R}
\overline{\Phi(\mathbf{u},b,\mathbf{a})}\psi(\mathbf{x}\cdot \mathbf{u}-b)K_{M^{*}}(\mathbf{a},\mathbf{x})
\,db d\va d\vu
}\\
=&
\int_{\mathbb R^n}
f(\mathbf{x}) d\vx
\int_{\mathbb{S}^{n-1}}
\int_{\mathbb R^n}
\int_{\mathbb R}
{\Phi(\mathbf{u},b,\mathbf{a})}
\overline{\psi(\mathbf{x}\cdot \mathbf{u}-b)}
\overline{K_{M^{*}}(\mathbf{a},\mathbf{x})}
\,db d\va d\vu
\end{align*}\begin{align*}
=&
\int_{\mathbb{S}^{n-1}}
\int_{\mathbb R^n}
\int_{\mathbb R}
{\Phi(\mathbf{u},b,\mathbf{a})}
\,db d\va d\vu
\int_{\mathbb R^n}
f(\mathbf{x})\overline{\psi(\mathbf{x} \cdot \mathbf{u}-b)}
{K_M(\mathbf{x},\mathbf{a})}
\,d\mathbf{x}
\\
=&
\int_{\mathbb{S}^{n-1}}
\int_{\mathbb R^n}
\int_{\mathbb R}
DS_\psi^M f(\mathbf{u},b,\mathbf{a}){\Phi(\mathbf{u},b,\mathbf{a})}
\,db d\va d\vu=\left\langle DS_\psi^M f(\vu,b,\va),\Phi(\vu,b,\va) \right\rangle_{\vu,b,\va}.
\end{align*}

\end{proof}

Analogously, we define the directional short-time special affine Fourier synthesis operator $(DS_{\psi}^{M})^{*}F$ of $F\in\mathcal S'(\mathbb Y^{2n})$ with respect to $\psi\in\mathcal S(\mathbb R)$ by
\begin{equation} \label{D2}
    \left\langle
\left(DS_{\psi}^{M}\right)^{*}F,\varphi
\right\rangle
:=
\left\langle
F,
\overline{
DS_{\psi}^{M}\overline{\varphi}
}
\right\rangle,
\quad
\varphi\in\mathcal S(\mathbb R^{n}).
\end{equation}
By Thm.~\ref{ContinuityDSTSAFT}, the distribution $(DS_{\psi}^{M})^*F$ is well defined via the duality \eqref{D2}. Furthermore, by (\cite{treves}, Corollary of Prop. 19.5), the mapping
$
\left(DS_{\psi}^{M}\right)^*:\mathcal{S}'(\mathbb{Y}^{2n})\longrightarrow \mathcal{S}'(\mathbb{R}^n) $
is continuous.

The following result extends the reconstruction formula \eqref{Inversionsynthesisnr} to the space of tempered distributions.

\begin{proposition}\label{relation1}
Let $\eta\in\mathcal S(\mathbb R)$ be a synthesis window for $\psi\in\mathcal S(\mathbb R) \setminus \{0\}$. Then, the Calderón reproducing formula holds
\begin{equation*}
    \frac{1}{(\eta,\psi)_{L^2(\mathbb R)}}
\left(
\left(DS_\eta^M\right)^*
\circ
DS_\psi^M
\right)
=
\operatorname{Id}_{\mathcal S'(\mathbb R^n)}.
\end{equation*}
\end{proposition}

\begin{proof}  Let $f\in\mathcal S'(\mathbb R^n)$. By Definitions \eqref{D1} and \eqref{D2}, and the inversion formula \eqref{Inversionsynthesisnr}, we obtain
\begin{align*}
    \left\langle
\left(DS_\eta^M\right)^*
DS_\psi^M f,\varphi
\right\rangle
&=
\left\langle
DS_\psi^M f,
\overline{DS_\eta^M\overline{\varphi}}
\right\rangle
=
\left\langle
f,
\overline{
\left(DS_\psi^M\right)^*DS_\eta^M\overline{\varphi}
}
\right\rangle
\\
&=
\left\langle
f,
\overline{{(\psi,\eta)}_{L^2(\mathbb R)}}\,\varphi
\right\rangle
=
{(\eta,\psi)}_{L^2(\mathbb R)}
\left\langle
f,\varphi
\right\rangle
\end{align*}
for all $\varphi\in\mathcal S(\mathbb R^n)$.
\end{proof}

The following result establishes an important relationship between the DSTSAFT and the SAFT.

{For each $(\mathbf{u},b)\in \mathbb{S}^{n-1}\times\mathbb{R}$, the function
$\mathbf{x}\mapsto \overline{\psi(\mathbf{x}\cdot \mathbf{u}-b)}$ is a multiplier on
$\mathcal{S}(\mathbb{R}^n)$; hence, for every
$f\in\mathcal{S}'(\mathbb{R}^n)$, the product
$f(\cdot)\overline{\psi((\cdot)\cdot \mathbf{u}-b)}$ is well defined as a tempered distribution.}

 \begin{proposition}\label{mixedSAFT}
Let $f \in \mathcal{S}'(\mathbb{R}^n)$ and $\psi \in \mathcal{S}(\mathbb{R})$. Then, for $\Phi \in \mathcal{S}(\mathbb{Y}^{2n})$,
\begin{equation} \label{CN}
\left\langle DS_\psi^M f, \Phi \right\rangle
= \int_{\mathbb{S}^{n-1}} \int_{\mathbb{R}}
\left\langle S_M\left(f(\cdot)\,\overline{\psi((\cdot)\cdot \mathbf{u} - b)}\right)(\mathbf{a}),\ \Phi(\mathbf{u},b,\mathbf{a}) \right\rangle_{\mathbf{a}}
\, db\, d\mathbf{u}.
\end{equation}
Furthermore, $DS_\psi^M f \in C^\infty(\mathbb{S}^{n-1}\times\mathbb{R}, \mathcal{S}'(\mathbb{R}^n))$ and  for every $\varphi\in\mathcal S(\mathbb R^n)$ the scalar-valued function
\begin{equation*}
        (\mathbf u,b)\mapsto
\left\langle DS_\psi^M f(\mathbf u,b,\cdot),\varphi\right\rangle
\end{equation*} is of slow growth on $\mathbb{S}^{n-1}\times\mathbb{R}$.
\end{proposition}
\begin{proof}
Definitions \eqref{D1} and \eqref{synthesis}, together with Fubini's theorem, yield 
\begin{align*}
    \left\langle DS_\psi^M f, \Phi\right\rangle
&=
\left\langle
f(\mathbf x),
\overline{\left(DS_\psi^M\right)^*\overline{\Phi}(\mathbf x)}
\right\rangle_{\vx}
\\
&=
\left\langle
f(\mathbf x),
\overline{
\int_{\mathbb{S}^{n-1}}
\int_{\mathbb R^n}
\int_{\mathbb R}
\overline{\Phi(\mathbf{u},b,\mathbf{a})}
\psi(\mathbf{u}\cdot\mathbf{x}-b)
K_{M^{*}}(\mathbf{a},\mathbf{x})
\,db d\va d\vu}
\right\rangle_{\vx}
\\
&=
\left\langle
f(\mathbf{x}),
\int_{\mathbb{S}^{n-1}}
\int_{\mathbb R}
\overline{\psi(\mathbf x\cdot\mathbf u-b)} db\,d\mathbf u
\int_{\mathbb R^n}
K_M(\mathbf x,\mathbf a)\,
\Phi(\mathbf u,b,\mathbf a)\,
d\mathbf a
\right\rangle_{\vx}\\
&= \left\langle f(\mathbf x),\ \int_{\mathbb S^{n-1}}\int_{\mathbb R}
\overline{\psi(\mathbf x \cdot \mathbf{u} -b) \, S_{M^*} \left(\overline {\Phi}(\mathbf{u},b,\cdot)\right)(\mathbf x)} \, db\,d\mathbf{u} \right\rangle_{\vx},
\end{align*}
where $M^{*}$ is given by \eqref{M*}.

According to the structural theorem (\cite{Schwartz}, page 239, Thm. VI) for $f \in \mathcal{S}'(\mathbb{R}^n)$, there exists a continuous function $h$ with at most polynomial growth and a multi-index $\boldsymbol{\beta} \in \mathbb N_0^n$ such that $f = \partial^{\boldsymbol{\beta}} h$ in the distributional sense. {The polynomial growth of $h$ and the rapid decay of the remaining factors justify the use of Fubini's theorem. Hence, by}  substituting {$f=\partial^\beta h$} into the previous relation, we obtain
\begin{align*}
\left\langle DS_\psi^M f, \Phi \right\rangle
&= (-1)^{|\boldsymbol{\beta}|} \int_{\mathbb R^n} h(\mathbf x)d\vx \frac{\partial^{|\boldsymbol{\beta}|}}{\partial \mathbf x^{\boldsymbol{\beta}}}\left( \int_{\mathbb S^{n-1}}\int_{\mathbb R} \overline{\psi(\mathbf{x}\cdot \mathbf{u}-b)\, S_{M^*}\left(\overline{\Phi}(\mathbf{u},b,\cdot)\right)(\mathbf x)}\right)\, db\,d\mathbf{u}  \\
&= (-1)^{|\boldsymbol{\beta}|} \int_{\mathbb S^{n-1}}\int_{\mathbb R} db\,d\mathbf{u} \int_{\mathbb R^n} h(\mathbf x)\, \frac{\partial^{|\boldsymbol{\beta}|}}{\partial \mathbf x^{\boldsymbol{\beta}}} \left(\overline{\psi (\mathbf x\cdot \mathbf{u} - b)\, S_{M^*} \left(\overline{ \Phi}(\mathbf{u},b,\cdot)\right)(\mathbf{x} )}  \right)\, d\mathbf x \\
&= \int_{\mathbb S^{n-1}}\int_{\mathbb R} \left\langle f(\mathbf x),\ \overline{\psi(\mathbf x\cdot \mathbf{u}-b)\, S_{M^*}\left(\overline{\Phi}(\mathbf{u},b,\cdot)\right)(\vx)} \right\rangle_{\mathbf x}\, db\,d\mathbf{u} \\
&= \int_{\mathbb S^{n-1}}\int_{\mathbb R} \left\langle S_M\left(f(\cdot)\,\overline{\psi((\cdot)\cdot \mathbf{u} - b)}\right)(\mathbf a),\ \Phi(\mathbf{u},b,\mathbf a) \right\rangle_{\mathbf a}\, db\,d\mathbf{u},
\end{align*}
which proves the relation \eqref{CN}. We now prove that, for every fixed $\varphi \in \mathcal S(\mathbb R^n)$, the function
\begin{equation} \label{FH}
\begin{split}
  (\mathbf u,b) \longmapsto \left\langle S_M\left(f(\cdot)\,\overline{\psi((\cdot)\cdot \mathbf u - b)}\right)(\mathbf a),\ \varphi(\mathbf a)\right\rangle_{\mathbf a}=\left\langle f(\mathbf x),\ \overline{\psi(\mathbf x\cdot \mathbf{u}-b)\,S_{M^*}\overline{\varphi}(\mathbf x)} \right\rangle_{\mathbf x},  
\end{split} 
\end{equation}
where $(\mathbf u,b) \in \mathbb S^{n-1}\times \mathbb R$ is smooth and has at most polynomial growth on $\mathbb S^{n-1}\times \mathbb R$. 

Since $\psi \in \mathcal{S}(\R)$ and $S_{M^*}\overline{\varphi} \in \mathcal{S}(\R^n)$, the derivatives of $\overline{\psi(\mathbf x\cdot \mathbf{u}-b)\, S_{M^*}\overline{\varphi}(\mathbf x)}$ with respect to the variables $\vu$ and $b$  still belong to $\mathcal{S}(\R^{n})$ as functions of $\vx$. By (\cite{treves}, Thm. 27.1), we may conclude that the function defined by \eqref{FH} is smooth on $ \mathbb{S}^{n-1} \times \mathbb{R}$. Furthermore, since $f\in \mathcal{S}' (\mathbb{R}^n)$, there exist $\nu \in \mathbb{N}_0$ and a constant $C'>0$ such that,
\begin{equation} \label{18}
    \left | \langle f,\varphi \rangle \right | \leq C' \rho_{\nu}(\varphi),
\end{equation}
for all $\varphi \in \mathcal{S}(\mathbb{R}^n)$. Then, by \eqref{FH} and \eqref{18}, we have

\begin{align*}
&\left |\left\langle DS_\psi^M f(\mathbf u,b, \mathbf a),\ \varphi(\mathbf a) \right\rangle_{\mathbf a}\right | =\left |\left\langle f(\mathbf x),\ \overline{\psi(\mathbf x\cdot \mathbf u - b)\, S_{M^*}\overline{\varphi}(\vx)} \right\rangle_{\mathbf x}\right | \\
\leq&
C'
\rho_{\nu}
\left(
\overline{\psi(\mathbf x\cdot\mathbf u-b)\, S_{M^{*}}\overline{\varphi}(\vx)}
\right) =C' \sup_{\substack{\mathbf x\in\mathbb R^n\\|\boldsymbol\alpha|\leq\nu}}
\left(1+|\mathbf x|^2\right)^{\nu/2}
\left|
\partial_{\mathbf x}^{\boldsymbol\alpha}
\left(\overline{\psi(\mathbf x\cdot\mathbf u-b)\,S_{M^{*}}\overline{\varphi}(\vx)}
\right)
\right|\\
\leq&C \sum_{\boldsymbol\beta\leq\boldsymbol\alpha}
\binom{\boldsymbol\alpha}{\boldsymbol\beta} \sup_{\substack{\mathbf x\in\mathbb R^n\\|\boldsymbol\alpha|\leq\nu}}
\left(1+|\mathbf x|^2\right)^{\nu/2}
\left|\partial_{\mathbf x}^{\boldsymbol\beta}
\left(\overline{\psi(\mathbf x\cdot\mathbf u-b)}\right)\right|
\,\left|\partial_{\mathbf x}^{\boldsymbol\alpha-\boldsymbol\beta}
\left(\overline{{S_{M^{*}}\overline{\varphi}(\vx)}}\right)\right|\\
\lesssim&
\rho_{\nu}(\psi)\sum_{
\boldsymbol\beta\leq\boldsymbol\alpha}
\sup_{\substack{\mathbf x\in\mathbb R^n\\|\boldsymbol\alpha|\leq\nu}}
\left(1+|\mathbf x|^2\right)^{\nu/2}
\left|
\partial_{\mathbf x}^{\boldsymbol\alpha-\boldsymbol\beta}\left(
\overline{{S_{M^{*}}\overline{\varphi}(\vx)}}\right)\right| \lesssim \rho_{\nu}(\psi) \rho_{\nu}\left(\overline{{S_{M^{*}}\overline{\varphi}}}\right).  
\end{align*}

Then, for every fixed $\varphi \in \mathcal S(\mathbb R^n)$, by means of the standard identification, the function given by \eqref{FH} can be viewed as an element of $\mathcal{S}'(\mathbb{S}^{n-1}\times \R)$ via the action 
\begin{align*}
&\left\langle \left\langle S_M\left(f(\cdot)\,\overline{\psi((\cdot)\cdot \mathbf u - b)}\right)(\mathbf a),\ \varphi(\mathbf a)\right\rangle_{\mathbf a},\ \Psi(\mathbf u,b) \right\rangle_{\mathbf u,b}\\
=& \int_{\mathbb S^{n-1}}\int_{\mathbb R} \left\langle S_M\left(f(\cdot)\,\overline{\psi((\cdot)\cdot \mathbf u - b)}\right)(\mathbf a),\ \varphi(\mathbf a) \right\rangle_{\mathbf a}\, \Psi(\mathbf u,b)\, db\, d\mathbf u, \quad \Psi \in \mathcal{S}(\mathbb{S}^{n-1}\times \R).
\end{align*}

The previous relation, together with \eqref{CN}, yields, for every $\varphi\in\mathcal{S}(\R^n)$ and $\Psi \in \mathcal{S}(\mathbb{S}^{n-1}\times \R)$,
\begin{align*}
&\left\langle \left \langle DS_\psi^M f, \varphi\right\rangle,\Psi \right\rangle=\left\langle \left\langle S_M\left(f(\cdot)\,\overline{\psi((\cdot)\cdot \mathbf u - b)}\right)(\mathbf a),\ \varphi(\mathbf a)\right\rangle_{\mathbf a},\ \Psi(\mathbf u,b) \right\rangle_{\mathbf u,b}
\end{align*}
from which we obtain
\begin{align*}
\left\langle DS_\psi^M f(\mathbf u,b, \mathbf a),\ \varphi(\mathbf a) \right\rangle_{\mathbf a}
&= \left\langle S_M\left(f(\cdot)\,\overline{\psi((\cdot)\cdot \mathbf u - b)}\right)(\mathbf a),\ \varphi(\mathbf a) \right\rangle_{\mathbf a} \\
&= \left\langle f(\mathbf x),\ \overline{\psi(\mathbf x\cdot \mathbf u - b)\, S_{M^*}\overline{\varphi}(\vx)}\right\rangle_{\mathbf x}.
\end{align*}
Thus, $DS_\psi^M f \in C^\infty\left(\mathbb S^{n-1}\times \mathbb R,\ \mathcal S'(\mathbb R^n)\right)$ and is of at most polynomial growth on $\mathbb S^{n-1}\times \mathbb R$.
\end{proof}

As a consequence of Prop. \ref{mixedSAFT} and the Calderón reproducing formula, we obtain the following desingularization formula.

\begin{corollary}[Desingularization formula]\label{cor:desing-SAFT}

Let $M = {(A,B,\mathbf{p};C,D,\mathbf{q})}$, $\psi \in \mathcal{S}(\R)\setminus \{0\}$, and $f \in \mathcal S'(\mathbb R^n)$. Then, for every $\varphi \in \mathcal S(\mathbb R^n)$, we have
\begin{equation*}
\langle f,\varphi\rangle
= \frac{1}{(\eta,\psi)_{L^2(\mathbb{R})}} \int_{\mathbb S^{n-1}}\int_{\mathbb R}
\left\langle S_M\left(f(\cdot)\,\overline{\psi((\cdot)\cdot \mathbf u-b)}\right)(\mathbf a),\ \overline{DS_\eta^M\overline{\varphi}(\mathbf u,b,\mathbf a)}\right\rangle_{\mathbf a}\, db\, d\mathbf u,
\end{equation*}
where $\eta \in \mathcal S(\mathbb R)$ is a synthesis window for $\psi$.
\end{corollary}

\begin{proof} Combining Propositions \ref{relation1} and \ref{mixedSAFT}, for $\varphi \in \mathcal S(\mathbb R^n)$, we obtain

\begin{align*}
    \langle f,\varphi\rangle
&=
\frac{1}{(\eta,\psi)_{L^2(\mathbb R)}}
\left\langle
\left(DS_\eta^M\right)^*DS_\psi^M f,\varphi
\right\rangle=
\frac{1}{(\eta,\psi)_{L^2(\mathbb R)}}
\left\langle
DS_\psi^M f,\overline{DS_\eta^M\overline{\varphi}}
\right\rangle
\\
&= \frac{1}{(\eta,\psi)_{L^2(\mathbb{R})}} \int_{\mathbb S^{n-1}}\int_{\mathbb R}
\left\langle S_M\left(f(\cdot)\,\overline{\psi((\cdot)\cdot \mathbf u-b)}\right)(\mathbf a),\ \overline{DS_\eta^M\overline{\varphi}(\mathbf u,b,\mathbf a)}\right\rangle_{\mathbf a}\, db\, d\mathbf u. 
\end{align*}
\end{proof}


\section{Proof of main results} \label{LS}

This section is devoted to the proofs of Theorems \ref{ContinuityDSTSAFT} and \ref{ContinuityDSTSAFTsynthesis}. The arguments follow the same approach as those in \cite{ferizi1} and \cite{SanevaAtanasova2016}, where the corresponding continuity results were established for the DSTFRFT and the DSTFT, respectively. {The main additional feature in the present setting is the more general SAFT kernel $K_M$, whose affine and chirp terms have to be taken into account in the corresponding estimates.}

\subsection{Proof of Theorem \ref{ContinuityDSTSAFT}}

It suffices to show that for arbitrary  $s, r, m, k \in \mathbb{N}_0, \vl \in \mathbb{N}_0^n$, there exist $\nu,\tau \in \mathbb{N}_0$ and a constant $C>0$ such that for all {$f\in \mathcal S(\mathbb R^n)$} and ${\psi\in \mathcal S(\mathbb R)}$
\begin{equation*}
    \rho_{s,r}^{\vl,m,k}\left(DS_\psi^M f\right) \leq C\rho_\nu(f)\rho_\tau(\psi).
\end{equation*}

\noindent \textbf{Step 1.} {Using Leibniz's formula, we obtain}
\begin{align*}
&\left|\frac{\partial^{\vl}}{\partial \mathbf{a}^{\vl}} \frac{\partial^{m}}{\partial b^m} \left(DS_\psi^Mf(\mathbf{u},b,\mathbf{a})\right)\right|=
\left|\frac{\partial^{\vl}}{\partial \mathbf{a}^{\vl}} \frac{\partial^{m}}{\partial b^m}  \int_{\mathbb R^n}f(\mathbf{x})\overline{\psi(\mathbf{u}\cdot \mathbf{x}- b)} K_M(\mathbf{x},\mathbf{a})\,d\mathbf{x} \right|
\\
=&
\left| \frac{\partial^{\vl}}{\partial \mathbf{a}^{\vl}}  \int_{\mathbb R^n}f(\mathbf{x}) \overline{\psi^{(m)}(\mathbf{u}\cdot \mathbf{x}-b)}K_M(\mathbf{x},\mathbf{a}) \,d\mathbf{x}   \right|
\\
=&
\left| \int_{\mathbb R^n}f(\mathbf{x}) \overline{\psi^{(m)}(\mathbf{u}\cdot \mathbf{x}-b)} \sum_{|\boldsymbol{\alpha}|+|\boldsymbol{\beta}|\leq |\vl|}
c_{\boldsymbol{\alpha},\boldsymbol{\beta}}\,
\mathbf{a}^{\boldsymbol{\alpha}} \mathbf{x}^{\boldsymbol{\beta}}K_M(\mathbf{x},\mathbf{a}) d\mathbf{x} \right|
\\
\leq& 
\sum_{|\boldsymbol{\alpha}|+|\boldsymbol{\beta}|\leq |\vl|} \left| c_{\boldsymbol{\alpha},\boldsymbol{\beta}} \mathbf{a}^{\boldsymbol \alpha} \int_{\mathbb R^n} \mathbf{x}^{\boldsymbol{\beta}}f(\mathbf{x}) \overline{\psi^{(m)}(\mathbf{u}\cdot \mathbf{x}-b)} \, K_M(\mathbf{x},\mathbf{a})d\mathbf{x} \right|
\\
\leq& 
\left(1+|\mathbf{a}|^2\right)^{|\vl|/2}\sum_{|\boldsymbol{\alpha}|+|\boldsymbol{\beta}|\leq |\vl|} |c_{\boldsymbol{\alpha},\boldsymbol{\beta}}| \left| DS_{\psi^{(m)}}^M f_{\boldsymbol{\beta}}(\mathbf{u},b,\mathbf{a})\right|,
\end{align*}
where $f_{\boldsymbol{\beta}}(\mathbf{x})=\mathbf{x}^{\boldsymbol{\beta}} f(\mathbf{x}) {\in \mathcal S(\mathbb R^n)}$, and $c_{\boldsymbol{\alpha}, \boldsymbol{\beta}}$ are constants that depend only on $\vl$ {and the fixed
parameter $M$}. {Since multiplication by polynomials and differentiation preserve the Schwartz space, and the additional factor
$\left(1+|\mathbf a|^2\right)^{|\vl|/2}$ can be absorbed by replacing $s$ with
$s+|\vl|$, it suffices to consider the case}  $m=0$ and $\vl=(0,0,...,0) $.

\noindent \textbf{Step 2.} For a given $(b,\mathbf a)\in \mathbb R\times\mathbb R^n$, we denote 
\begin{equation*}
    F_{b,\va}: \mathbf y\mapsto DS_{\psi}^{M}f\left(\frac{\mathbf{y}}{|\mathbf{y}|},b,\mathbf a\right)= 
\int_{\mathbb{R}^n}
f(\mathbf{x})\, \overline{\psi\left(\frac{\mathbf{y} \cdot \vx}{|\mathbf{y}|} - b\right)}\,
K_M(\mathbf{x},\mathbf{a})\, d\mathbf{x}
\end{equation*} for $\mathbf{y} \in \R^{n}\setminus\{0\}$. {Since $F_{b,\mathbf a}$ is homogeneous of degree zero, the polar
coordinate representation of the Euclidean Laplacian gives,} for $\mathbf u\in\mathbb S^{n-1}$,
\begin{equation*}
    \Delta_{\mathbb{S}^{n-1}}D\mathcal S_{\psi}^{M}f(\mathbf u,b,\mathbf a)=\Delta F_{b,\va}(\mathbf u).
\end{equation*}

\noindent Since $\Delta_{\mathbb{S}^{n-1}}^{k}$ is a differential operator of order $2k$ on the unit sphere $\mathbb{S}^{n-1}$, it follows that
\begin{align*}
&\left|
\Delta_{S^{n-1}}^{k}
DS_{\psi}^{M}
f(\mathbf u,b,\mathbf a)
\right|
=  
\left|
\sum_{\substack{0\le j\le2k\\
|\beta|\le2k}}
P_{j,\beta,k}(\mathbf u) \int_{\mathbb{R}^n} \mathbf{x}^{\boldsymbol{\beta}} f(\mathbf{x}) \overline{\psi^{(j)}(\mathbf{x} \cdot \mathbf{u}-b)} K_{M}(\mathbf{x},\mathbf{a}) d\mathbf{x}
\right|
\\
\le&
\sum_{\substack{0\le j\le2k\\
|\beta|\le2k}}
|P_{j,\beta,k}(\mathbf u)|
\,
\left|
DS_{\psi^{(j)}}^{M}
f_{\beta}(\mathbf u,b,\mathbf a)
\right|
\lesssim
\sum_{\substack{0\le j\le2k\\
|\beta|\le2k}}
\left|
DS_{\psi^{(j)}}^{M}
f_{\beta}(\mathbf u,b,\mathbf a)
\right|,
\end{align*}
where $P_{j,\boldsymbol{\beta},k}(\mathbf{u})$ are polynomials depending only on $k$, $f_{\boldsymbol{\beta}}(\mathbf{x})=\mathbf{x}^{\boldsymbol{\beta}}f(\mathbf{x}) \in\mathcal S(\mathbb R^n)$, and $\psi^{(j)}\in\mathcal S(\mathbb R)$. {Since $\mathbb S^{n-1}$ is compact, the polynomials
$P_{j,\beta,k}$ are uniformly bounded on $\mathbb S^{n-1}$.} Thus, it suffices to consider the case $k=0$.

By replacing $r$ and $s$, if necessary, by the next larger even integers, we may assume without loss of generality that both are even, since the corresponding seminorms with higher indices dominate the original ones.

\noindent \textbf{Step 3.} 
By Prop. \ref{relationwithFT}, we obtain
\begin{align*}
&\left(1+|b|^2\right)^{r/2}\left|DS_\psi^M f(\mathbf{u},b,\mathbf{a})\right|
\\=&
(2\pi)^{-1/2} \left|
\int_{\mathbb R}
S_M\!\left(f(\cdot)e^{-i(\cdot)(\xi \mathbf{u})}\right)(\mathbf{a})
\overline{\widehat{\psi}(\xi)}
\left(
1-\frac{\partial^2}{\partial\xi^2}
\right)^{r/2}
\left(e^{i\xi b}\right)
\,d\xi \right|\\ =&
(2\pi)^{-1/2}
\left|\int_{\mathbb R}
\left( 
1-\frac{\partial^2}{\partial\xi^2}
\right)^{r/2}
\left(
S_M\!\left(f(\cdot)e^{-i(\cdot)(\xi\mathbf{u})}\right)(\mathbf{a})
\overline{\widehat{\psi}(\xi)}
\right)
e^{i\xi b}
\,d\xi \right|.\end{align*}
The last equality follows from repeated integration by parts, with all boundary terms vanishing as a consequence of the rapid decay of the Schwartz functions involved. Therefore, an application of Prop. \ref{relationwithFT} yields
\begin{align*}
&\left(1+|b|^2\right)^{r/2} \left|DS_\psi^M f(\mathbf{u},b,\mathbf{a}) \right|\\
\leq& {(2\pi)^{-1/2} \sum_{j\leq r/2} \binom{r/2}{j}\left| \int_{\mathbb R}
\frac{\partial^{2j}}{\partial\xi^{2j}}
\left(
S_M\!\left(f(\cdot)e^{-i(\cdot)(\xi\mathbf{u})}\right)(\mathbf{a})
\overline{\widehat{\psi}(\xi)}
\right)
e^{i\xi b}
\,d\xi\right|}\end{align*}\begin{align*}
\leq& {(2\pi)^{-1/2} \sum_{j\leq r/2}\sum_{d\leq 2j} \binom{r/2}{j} \binom{2j}{d}\left| \int_{\mathbb R}
\frac{\partial^{d}}{\partial\xi^{d}}
\left(
S_M\!\left(f(\cdot)e^{-i(\cdot)(\xi\mathbf{u})}\right)(\mathbf{a})\right) \frac{\partial^{2j-d}}{\partial\xi^{2j-d}}\left(
\overline{\widehat{\psi}(\xi)}\right)
e^{i\xi b}
\,d\xi\right|}\\
\lesssim & {(2\pi)^{-1/2} \sum_{j\leq r/2}\sum_{d\leq 2j}\left| \int_{\mathbb R}
S_M\!\left((\vx\cdot \vu)^df(\vx)e^{-i\vx \cdot (\xi\mathbf{u})}\right)(\mathbf{a})
\overline{\mathcal{F}(t^{2j-d} \psi(t))(\xi)}
e^{i\xi b}
\,d\xi\right|}\\
\lesssim & {(2\pi)^{-1/2} \sum_{j\leq r/2}\sum_{d\leq 2j} \sum_{|\boldsymbol{\beta}|\leq d} \left| \int_{\mathbb R}
S_M\!\left(\vx^{\boldsymbol{\beta}} f(\vx)e^{-i\vx \cdot (\xi\mathbf{u})}\right)(\mathbf{a})
\overline{\mathcal{F}(t^{2j-d} \psi(t))(\xi)}
e^{i\xi b}
\,d\xi\right|}\\
{\lesssim}&\sum_{j\leq r/2}\sum_{d\leq 2j} \sum_{|\boldsymbol{\beta}|\leq d}
\left|DS_{\psi_{2j-d}}^M f_{\boldsymbol{\beta}}(\mathbf{u},b,\mathbf{a})\right|,
\end{align*}
where  $f_{\boldsymbol{\beta}} {(\vx)}=\mathbf{x}^{\boldsymbol{\beta}}f(\mathbf{x}) \in \mathcal{S}(\mathbb{R}^n)$ and $\psi_{2j-d}(t)=t^{2j-d}\psi(t) {\in} \mathcal{S}(\mathbb{R})$. {Thus, we may reduce the proof to the case $r=0$.}

\noindent \textbf{Step 4.} 
To obtain the required decay with respect to the variable $\mathbf{a}$, we first estimate the polynomial weight in terms of $B^{-1}\mathbf{a}$. Since $\mathbf{a}=B \left(B^{-1}\mathbf{a} \right)$, the spectral norm inequality gives $|\mathbf{a}| \leq \|B \|_2 | B^{-1}\mathbf{a}|$. Consequently, 
\begin{align*}
    \left(1 + |\mathbf{a}|^2\right)^{s/2} \left|
DS_{\psi}^{M}f(\mathbf u,b,\mathbf a)
\right| \leq & \left(1 + \|B \|_2^2 \, \left| B^{-1}\mathbf{a}\right|^2\right)^{s/2} \left|
DS_{\psi}^{M}f(\mathbf u,b,\mathbf a)
\right|\\
\leq& \left(1 + \|B \|_2^2 \right)^{s/2}\, \left(1 + \left| B^{-1}\mathbf{a}\right|^2\right) ^{s/2}\left|
DS_{\psi}^{M}f(\mathbf u,b,\mathbf a)
\right|,
\end{align*}
where $\|\cdot\|_2$ is the spectral norm of a matrix. Since $s$ is even and
$\chi_{M}(\mathbf{x})=\frac{1}{2}\mathbf{x}^{T}B^{-1}A\mathbf{x}+\mathbf{p}^{T}B^{-1}\mathbf{x}$, repeated integration by parts yields
\begin{align*}
   & \left(1 + |\mathbf{a}|^2\right)^{s/2} \left|
DS_{\psi}^{M}f(\mathbf u,b,\mathbf a)
\right| \\
\lesssim &\left(1 + \left| B^{-1}\mathbf{a}\right|^2\right) ^{s/2}
\left|
\int_{\mathbb R^{n}}
f(\mathbf x)
\overline{\psi(\mathbf u\cdot\mathbf x-b)}
e^{i\chi_{M}(\mathbf x)}
e^{-i\mathbf x^{T}B^{-1}\mathbf a}
\,d\mathbf x
\right| \\
\lesssim&
\left|
\int_{\mathbb R^n}
f(\mathbf x)
\overline{\psi(\mathbf u\cdot\mathbf x-b)}
e^{i\chi_M(\mathbf x)}
\left(
1-\sum_{k=1}^n\frac{\partial^2}{\partial x_k^2}
\right)^{s/2}
\left(e^{-i\mathbf x^TB^{-1}\mathbf a}\right)
\,d\mathbf x
\right|
\\
=&
\left|
\int_{\mathbb R^n}
\left(
1-\sum_{k=1}^n\frac{\partial^2}{\partial x_k^2}
\right)^{s/2}
\left(
f(\mathbf x)
e^{i\chi_M(\mathbf x)}
\overline{\psi(\mathbf u\cdot\mathbf x-b)}
\right)
e^{-i\mathbf x^TB^{-1}\mathbf a}
\,d\mathbf x
\right|.
\end{align*}
Expanding the differential operator yields powers of the Laplacian up to
order $s/2$, and hence derivatives of total order at most $s$. Applying the Leibniz rule to these derivatives, and using the fact that derivatives of the chirp factor $e^{i \chi_{M}(\mathbf{x})}$ are polynomials in $\mathbf{x}$ multiplied by the same chirp factor, we obtain 

\begin{align*}
& \left(1 + |\mathbf{a}|^2\right)^{s/2} \left|
DS_{\psi}^{M}f(\mathbf u,b,\mathbf a)
\right|\\
\lesssim&
\sum_{|\boldsymbol\alpha|\leq s/2}
\left|
\int_{\mathbb R^n}
\partial_{\mathbf x}^{2\boldsymbol\alpha}
\left(
f(\mathbf x)
e^{i\chi_M(\mathbf x)}
\overline{\psi(\mathbf u\cdot\mathbf x-b)}
\right)
e^{-i\mathbf x^TB^{-1}\mathbf a}
\,d\mathbf x
\right|
\\
\lesssim&
\sum_{|\boldsymbol\alpha|\leq s/2}\, \sum_{|\boldsymbol\beta|, |\boldsymbol\mu|, j \leq 2 |\boldsymbol{\alpha}|}
\int_{\mathbb R^n}\left|
\partial^{\boldsymbol\beta}
f(\mathbf x) \vx^{\boldsymbol{\mu}}
\overline{\psi^{(j)}(\mathbf u\cdot\mathbf x-b)}
\right|
\,d\mathbf x\\
\le& \sum_{|\boldsymbol\alpha|\leq s/2}\, \sum_{|\boldsymbol\beta|, |\boldsymbol\mu|, j \leq 2 |\boldsymbol{\alpha}|}
\int_{\mathbb R^n}
(1+|\mathbf x|)^{|\boldsymbol\mu|+n+1}
(1+|\mathbf x|)^{-n-1}
\left|
\partial^{\boldsymbol\beta}f(\mathbf x)
\right|\left|\overline{\psi^{(j)}(\mathbf u\cdot\mathbf x-b)}\right|\,d\mathbf x
\\
\lesssim& \sum_{|\boldsymbol\alpha|\leq s/2}\, \sum_{|\boldsymbol\beta|, |\boldsymbol\mu|, j \leq 2 |\boldsymbol{\alpha}|}
\int_{\mathbb R^n}
(1+|\mathbf x|^2)^{(n+s+1)/2}
\left|\partial^{\boldsymbol\beta}f(\mathbf x)\right|(1+|\mathbf x|)^{-n-1}
\left|\overline{\psi^{(j)}(\mathbf u\cdot\mathbf x-b)}
\right|\,d\mathbf x
\\
\lesssim&\rho_{n+1+s}(f)\,
\rho_s(\psi) \int_{\mathbb R^n} (1+|\mathbf x|)^{-n-1}\,d\mathbf x \lesssim
\rho_{n+1+s}(f)\rho_{s}(\psi).
\end{align*}
Combining Steps~1--4 yields the required estimate and completes the proof.

\subsection{Proof of Theorem \ref{ContinuityDSTSAFTsynthesis}}

To prove continuity, it suffices to show that for every $\nu\in\mathbb{N}_0$,
there exist $\tau \in \mathbb{N}_0$, $N\in\mathbb{N}$, indices
$s_j,r_j,m_j,k_j\in\mathbb{N}_0$ and $\vl_j\in\mathbb{N}_0^n$,
$j=1,\ldots,N$,
and a constant $C>0$ such that
\begin{equation*}
    \rho_\nu\!\left((DS_\psi^M)^*\Phi\right)
\leq C\left(
\sum_{j=1}^{N}
\rho_{s_j,r_j}^{\,l_j,m_j,k_j}(\Phi)
\right)\rho_\tau(\psi),
\end{equation*}
for every $\Phi\in\mathcal{S}(\mathbb{Y}^{2n})$ and
$\psi\in\mathcal{S}(\mathbb{R})$.

One may verify that, for every $\Phi \in \mathcal{S}(\mathbb{Y}^{2n})$ and $ \psi\in \mathcal{S}(\mathbb{R})$, $\left(D S_{\psi}^{M}\right)^{*} \Phi$ is a smooth function. We may suppose that $\nu$ is even. From the definition of the synthesis operator and Leibniz's formula, it follows that, for every $|\boldsymbol{\alpha}|\leq \nu$,

\begin{align*}
&\left(1+|\mathbf x|^{2}\right)^{\nu/2}
\left|
\partial_{\mathbf x}^{\boldsymbol{\alpha}}(DS_{\psi}^{M})^{*}\Phi(\mathbf x) \right|\\=&\left(1+|\mathbf x|^{2}\right)^{\nu/2}
\left|
\partial_{\mathbf x}^{\boldsymbol{\alpha}}\left(\int_{\mathbb{S}^{n-1}}
\int_{\mathbb R^n}
\int_{\mathbb R}
\Phi(\mathbf{u},b, \mathbf{a})
\psi(\mathbf{u}\cdot \mathbf{x}-b)
K_{M^{*}}(\mathbf{a},\mathbf{x})
\,db d\va d\mathbf{u}\right) \right|\\
\lesssim &
\left(1+|\mathbf x|^{2}\right)^{\nu/2}
\sum_{|\boldsymbol{\beta}|,|\boldsymbol{\gamma}|,|\boldsymbol{\mu}|\leq |\boldsymbol{\alpha}|}
\left|
\int_{\mathbb S^{n-1}}
\int_{\mathbb R^n}
\int_{\mathbb R}
\Phi(\mathbf u,b,\mathbf a)
\mathbf x^{\boldsymbol{\beta}}
\mathbf a^{\boldsymbol{\mu}} \vu^{\gamma}
\psi^{(|\boldsymbol{\gamma}|)}
(\mathbf u\cdot\mathbf x-b)
K_{M^{*}}(\mathbf a,\mathbf x)
\,db d\va d\mathbf{u}
\right|
\\
\lesssim&
\sum_{|\boldsymbol{\beta}|,|\boldsymbol{\gamma}|,|\boldsymbol{\mu}|\leq |\boldsymbol{\alpha}|} \left(1+|\mathbf x|^{2}\right)^{\frac{\nu+2|\boldsymbol{\beta}|}{2}}
\Bigg|
\int_{\mathbb S^{n-1}}
\int_{\mathbb R^n}
\int_{\mathbb R}
\Phi(\mathbf u,b,\mathbf a)
\mathbf a^{\boldsymbol{\mu}} \vu^{\gamma}
\psi^{(|\boldsymbol{\gamma}|)}
(\mathbf u\cdot\mathbf x-b)\\&
e^{ 
-\frac{i}{2} \left( 
\mathbf{a}^T D B^{-1} \mathbf{a}
+ \boldsymbol{\Omega} B^{-T} \mathbf{a}
- 2 \mathbf{a}^T B^{-T} \mathbf{x}
\right) 
}
\,db d\va d\mathbf{u}
\Bigg|
\end{align*}\begin{align*}
\leq& \sum_{|\boldsymbol{\beta}|,|\boldsymbol{\gamma}|,|\boldsymbol{\mu}|\leq |\boldsymbol{\alpha}|} \left(1+\|B\|_{2}^{2}\right)^{\frac{\nu+2|\boldsymbol{\beta}|}{2}} \left(1+|B^{-T}\mathbf x|^{2}\right)^{\frac{\nu+2|\boldsymbol{\beta}|}{2}}
\Bigg|
\int_{\mathbb S^{n-1}}
\int_{\mathbb R^n}
\int_{\mathbb R}
\Phi(\mathbf u,b,\mathbf a)
\mathbf a^{\boldsymbol{\mu}}
\psi^{(|\boldsymbol{\gamma}|)}
(\mathbf u\cdot\mathbf x-b) \\&
e^{ 
-\frac{i}{2} \left( 
\mathbf{a}^T D B^{-1} \mathbf{a}
+ \boldsymbol{\Omega} B^{-T} \mathbf{a}
- 2 \mathbf{a}^T B^{-T} \mathbf{x}
\right) 
}
\,db d\va d\mathbf{u}
\Bigg| \\
\lesssim &\sum_{|\boldsymbol{\beta}|,|\boldsymbol{\gamma}|,|\boldsymbol{\mu}|\leq |\boldsymbol{\alpha}|} 
\Bigg|
\int_{\mathbb S^{n-1}}
\int_{\mathbb R^n}
\int_{\mathbb R}
\Phi(\mathbf u,b,\mathbf a)
\mathbf a^{\boldsymbol{\mu}} \psi^{(|\boldsymbol{\gamma}|)}
(\mathbf u\cdot\mathbf x-b) e^{-\frac{i}{2} \left(\mathbf{a}^T D B^{-1} \mathbf{a}
+ \boldsymbol{\Omega} B^{-T} \mathbf{a}\right)}\\& \left(
1-\sum_{k=1}^{n}
\frac{\partial^2}{\partial a_k^2}
\right)^{\frac{\nu+2|\boldsymbol{\beta}|}{2}}\left(
e^{i \mathbf{a}^T B^{-T} \mathbf{x}}\right)
\,db d\va d\mathbf{u}
\Bigg| \\
\lesssim &\sum_{|\boldsymbol{\beta}|,|\boldsymbol{\gamma}|,|\boldsymbol{\mu}|\leq |\boldsymbol{\alpha}|} 
\Bigg|
\int_{\mathbb S^{n-1}}
\int_{\mathbb R^n}
\int_{\mathbb R}
\left(
1-\sum_{k=1}^{n}
\frac{\partial^2}{\partial a_k^2}
\right)^{\frac{\nu+2|\boldsymbol{\beta}|}{2}}\left(\Phi(\mathbf u,b,\mathbf a)
\mathbf a^{\boldsymbol{\mu}}  e^{-\frac{i}{2} \left(\mathbf{a}^T D B^{-1} \mathbf{a}
+ \boldsymbol{\Omega} B^{-T} \mathbf{a}\right)}\right)\\& \psi^{(|\boldsymbol{\gamma}|)}
(\mathbf u\cdot\mathbf x-b) 
e^{i \mathbf{a}^T B^{-T} \mathbf{x}}
\,db d\va d\mathbf{u}
\Bigg| \\ \lesssim&\sum_{|\boldsymbol{\beta}|,|\boldsymbol{\gamma}|,|\boldsymbol{\mu}|\leq |\boldsymbol{\alpha}|} \, \sum_{|\mathbf{d}|\leq \nu + 2|\boldsymbol{\beta}|}
\Bigg|
\int_{\mathbb S^{n-1}}
\int_{\mathbb R^n}
\int_{\mathbb R} \partial_{\va}^{\mathbf{d}}
\Phi(\mathbf u,b,\mathbf a) \, 
Q_{\mathbf{d},\boldsymbol{\mu}}(\va)  e^{-\frac{i}{2} \left(\mathbf{a}^T D B^{-1} \mathbf{a}
+ \boldsymbol{\Omega} B^{-T} \mathbf{a}\right)}\\& \psi^{(|\boldsymbol{\gamma}|)}
\left(\mathbf u\cdot\mathbf x-b\right) 
e^{i \mathbf{a}^T B^{-T} \mathbf{x}}
\,db d\va d\mathbf{u}
\Bigg| \\\leq& \rho_{\nu}(\psi)\sum_{|\boldsymbol{\beta}|,|\boldsymbol{\gamma}|,|\boldsymbol{\mu}|\leq |\boldsymbol{\alpha}|} \, \sum_{|\mathbf{d}|\leq \nu + 2|\boldsymbol{\beta}|}
\int_{\mathbb S^{n-1}}
\int_{\mathbb R^n}
\int_{\mathbb R} \left|\partial_{\va}^{\mathbf{d}}
\Phi(\mathbf u,b,\mathbf a)\right| \, 
\left|Q_{\mathbf{d},\boldsymbol{\mu}}(\va) \right|   
\,db d\va d\mathbf{u} \\\lesssim& \rho_{\nu}(\psi) \sum_{|\boldsymbol{\beta}|,|\boldsymbol{\gamma}|,|\boldsymbol{\mu}|\leq |\boldsymbol{\alpha}|} \, \sum_{|\mathbf{d}|\leq \nu + 2|\boldsymbol{\beta}|} \rho^{\mathbf{d},0,0}_{s,2}(\Phi) 
\end{align*}

\noindent for some polynomials $Q_{\mathbf{d}, \boldsymbol{\mu}}(\mathbf{a})$ and $s \in \mathbb{N}_0$. This completes the proof.


\subsection*{Acknowledgments}

The work of Katerina Hadzi-Velkova Saneva was supported by the Ministry of Education and Science of North Macedonia under Grant No. 15-6171/25.



\end{document}